\documentclass[11pt,hyp,]{nyjm}
\usepackage{hyperref}
\hypersetup{nesting=true,debug=true,naturalnames=true}
\usepackage{graphicx,amssymb,upref}

\let\<\langle
\let\>\rangle
\newcommand{\doi}[1]{doi:\,\href{http://dx.doi.org/#1}{#1}}

\let\uml\"

\theoremstyle{definition}

\theoremstyle{remark}

\newtheorem{Theorem}{Theorem}
\newtheorem{Lemma}[Theorem]{Lemma}
\newtheorem{Corollary}[Theorem]{Corollary}

\newtheorem{Definition}{Definition}
\newtheorem{Example}{Example}

\newcommand{\Vol}{\operatorname{Vol}}

\newcommand{\cout}[1]{}
\def\co{\colon\thinspace}
\def\Fs{{\mathcal F}}

\def\Sl{\widetilde{\rm SL}_{2}}

\numberwithin{equation}{section}

\begin{document}


\title[Collapsing, growth \& Einstein metrics on smooth 4-manifolds]{Collapsing and group growth as obstructions to Einstein metrics \\ on some smooth 4-manifolds}

\author[H. Contreras Peruyero]{H. Contreras Peruyero}
\address{Centro de Ciencias Matem\'aticas, Universidad Nacional Aut\'onoma de M\'exico UNAM, Antigua Carretera a P\'atzcuaro \# 8701, Col. Ex Hacienda San Jos\'e de la Huerta, 58089, Morelia, Michoac\'an, Mexico}
\email{haydeeperuyero@im.unam.mx}

\author[P. Su\'arez-Serrato]{P. Su\'arez-Serrato}
\address{Instituto de Matem\'aticas, Universidad Nacional Aut\'onoma de M\'exico UNAM, Área de la Investigación Científica, Circuito exterior,
Ciudad Universitaria, 04510, Coyoac\'an, Mexico City, Mexico}
\email{pablo@im.unam.mx}  
\thanks{PSS thankfully acknowledges partial support by DGAPA-UNAM PAPIIT grant IN104819.}



\date{\today}


\commby{Prof. Claude LeBrun}

\begin{abstract}
We show that a combination of collapsing and excessive growth from the fundamental group impedes the existence of Einstein metrics on several families of smooth $4$--manifolds. These include infrasolvmanifolds whose fundamental group is not virtually nilpotent, most elliptic surfaces of zero Euler characteristic, geometrizable manifolds with hyperbolic factor geometries in their geometric decomposition, and higher graph $4$--manifolds without purely negatively curved pieces.  
\end{abstract}

\maketitle


\section{Introduction}

A Riemannian metric on a smooth manifold is Einstein if the Ricci tensor is a scalar multiple of the metric.  The Einstein equation for a Riemannian manifold ${\rm Ric}_g= D \cdot g$ is a nonlinear second-order system of partial differential equations.  Existence and obstruction results about solutions are in general hard to obtain. A smooth manifold is called Einstein if it admits a Riemannian metric $g$ solving the Einstein equation. In dimensions 2 and 3 a manifold is Einstein if and only if it has constant sectional curvature. 

A closed Einstein $4$--manifold $M$ must satisfy the Hitchin-Thorpe inequality \cite{hit,thor}:
\begin{eqnarray}\label{ht-int}
2\chi(M) \geq 3|\tau(M)|
\end{eqnarray}
Here $\chi(M)$ is the Euler characteristic and $\tau(M)$ the signature of $M$. 

W.\,P. Thurston defined a \emph{model geometry} as a complete simply connected Riemannian manifold $M$ such that the group of isometries acts transitively on $M$ and contains a discrete subgroup with a finite volume quotient \cite{Thu0}. These \emph{Thurston geometries} are fundamental for the geometrization of $3$--manifolds. In dimension four, Filipkiewicz classified the maximal geometries  \cite{Fi}. 

 The subtleties of the existence of Einstein metrics in dimension four is treated by various authors (see, for example Kotschick \cite{kot(1998),kot(1998-1)}, LeBrun \cite{leb-1996,leb-2001}, and Sambusetti \cite{sam}).
 
Our first result shows that some of the manifolds modeled on Thurston geometries in dimension four do not admit Einstein metrics:

\begin{Theorem}\label{thm::noEinst-Tgeom}
 Let $M$ be a closed orientable smooth $4$--manifold modeled on one of the Thurston geometries $ \mathbb{S}^{3}\times \mathbb{E},\,
 \mathbb{N}il^{3}\times \mathbb{E},\,
\mathbb{N}il^{4}, \, \mathbb{H}^{3}\times \mathbb{E},\, \Sl\times \mathbb{E},\, \mathbb{H}^{2}\times\mathbb{E}^{2},\,  \mathbb{S}ol^{4}_{1},\, \mathbb{S}ol^{4}_{0},$ or $\mathbb{S}ol^{4}_{m,n}.$ Then $M$ does not admit an Einstein metric.
\end{Theorem}

For the sake of brevity, we recommend interested readers consult Hillman's book for the definitions and examples of manifolds modeled on these Thurston geometries \cite{Hill}.
Recall that Jensen classified homogeneous and Einstein simply connected $4$--manifolds  \cite{Jensen}. Jablonski showed that a subgroup of the isometries of an Einstein solvmanifold acting properly discontinuously such that the homogeneous quotient manifold defined by identifying its orbits is a trivial group \cite{Jab}. 

Here, an {\it infrasolvmanifold} is a quotient $M=S / \Gamma $, with $S$ a $1$--connected solvable subgroup and $\Gamma$ is a closed torsion-free subgroup of ${\rm Aff}(S)$ such that  the component of the identity of $\Gamma$ lies in the maximal connected nilpotent normal subgroup of $S$ (its nilradical), $\Gamma / \Gamma \cap S$ has compact closure in ${\rm Aut}(S)$, and $M$ is compact. Moreover, Hillman proved that closed $4$--dimensional infrasolvmanifolds are diffeomorphic to geometric $4$--manifolds of solvable Lie type \cite{Hill2}. 

Theorem \ref{thm::noEinst-Tgeom}, therefore, yields a straightforward, yet noteworthy, proof that some infrasolvmanifolds of dimension four do not admit Einstein metrics:

\begin{Corollary}\label{cor::infrasolv-noEins}
Compact $4$--dimensional infrasolvmanifolds whose fundamental group is not virtually nilpotent do not admit Einstein metrics.
\end{Corollary}

Smooth, compact complex surfaces are another family of $4$--manifolds of interest. A complex surface $S$ is an elliptic surface if there is a holomorphic map $S\to C$ to a complex curve $C$ with general fiber an elliptic curve. We obtain a new proof of the following result for elliptic surfaces of zero Euler characteristic, found in the classic book by Besse \cite{Besse}:

\begin{Theorem}\label{thm::noEinst-EllSfaces}
Compact complex elliptic surfaces of zero Euler characteristic and which are not finitely covered by a $4$--torus or a $K3$ surface do not admit Einstein metrics. 
\end{Theorem}

This last result also recovers some results of LeBrun where he classified Einstein $4$--manifolds that admit a complex structure \cite[Lemma 5]{leb-2009}. The previous theorems can also be restated (for some geometric manifolds and elliptic surfaces) in terms of circle foliations and Seifert fibrations, as we shall see below. 

\begin{Definition} {\cite{Hill}}
A smooth oriented $4$--manifold $M$ is Seifert fibred if it is the total space of an orbifold bundle $M\to B$ with general fiber a torus over a $2$--orbifold $B$.
\end{Definition}

Under this definition, Seifert $4$--manifolds account for every compact complex surface diffeomorphic to some elliptic surface with $c_2=0$, and it also encompasses examples that do not admit any complex structure whatsoever \cite{Wall}. 

\begin{Theorem}\label{thm::noEinst-S1fols-Seifert4mfds}
Let $M$ be a closed orientable smooth $4$--manifold that is either foliated by geodesic circles or that admits a Seifert fibration. Assume $M$ is not finitely covered by $T^4$ of a $K3$ surface. Then $M$ does not admit an Einstein metric. 
\end{Theorem}

In the following results, we will cover several ways of generalizing the notions of Seifert fibration and of {\it graph manifold} (that have previously appeared in the work of several authors) from dimension three to dimension four. 

We will say that a manifold $M$ admits \emph{geometric decomposition} if there exists a finite collection of $2$--sided hypersurfaces $S$ such that each component of $M-S$ admits a Thurston geometry. 

Consider the following example of a manifold with a geometric decomposition into pieces with geometric structures modeled on  $\mathbb{H}^{2}\times\mathbb{E}^{2}$ and $\mathbb{H}^{3}\times \mathbb{E}$. 

\begin{Example}\label{ex::geoms}
Let $\Sigma$ be an orientable surface of genus$\geq 2$ and finite area with a single cusp and realize $\Gamma_1 \cong  \pi_1(\Sigma, \ast)$ as a Fuchsian group of isometries of $\mathbb{H}^{2}$. Truncate the cusp of $\mathbb{H}^{2} / \Gamma_1$ to obtain a compact orientable surface $M_1$ with a single boundary component $\partial M_1 \cong S^1$. Define:
\[
Z_1 := M_1 \times T^2
\]
Let $Y$ be an orientable, finite volume, hyperbolic $3$--manifold with a single cusp and realize $\Gamma_2 \cong \pi_1(Y, \ast)$ as a Kleinian group of isometries of $\mathbb{H}^{3}$. Truncate the cusp of $\mathbb{H}^{2} / \Gamma_2$ to obtain a compact orientable $3$--manifold with $M_2$ with a single boundary component $\partial M_2 \cong T^2$.
Define:
\[
Z_2 := M_2 \times S^1
\]

Now construct the manifold $M$ using the identity to glue $Z_1$ and $Z_2$ along their boundary $T^3$. This manifold $M$ has a geometric decomposition into the Thurston geometries $\mathbb{H}^{2}\times\mathbb{E}^{2}$ and $\mathbb{H}^{3}\times \mathbb{E}$. 
\end{Example}
 
The following result tells us that manifolds like the one constructed in Example \ref{ex::geoms} do not admit Einstein metrics.
 
 \begin{Theorem}\label{thm::noEinst-geoms}
 Let $M$ be a closed orientable smooth $4$--manifold that admits a geometric decomposition in the sense of Thurston into at least two pieces. If all of the geometric pieces of $M$ are modeled on
geometries $\mathbb{H}^{3}\times \mathbb{E},\, \Sl\times \mathbb{E},$ or $\, \mathbb{H}^{2}\times\mathbb{E}^{2},$ then $M$ does not admit an Einstein metric.
\end{Theorem}

The notion of \emph{graph manifold} has been successful in dimension 3, helping to distinguish different classes of geometric interest. Collapsing manifolds that are geometrizable in the sense of Thurston are one analog of this notion. There are also other possibilities, namely, decomposing a manifold into pieces with infranil fibrations and generalizing the notion of Seifert fibration. Consider the following definition, which extends ideas put forth by Frigerio-Lafont-Sisto \cite{FLS15}:

\begin{Definition}[\cite{CS19}]\label{def::HGM}
A compact, smooth $n$--manifold $M$, $n \geq 3$, is called a {\em higher graph manifold} if it is constructed as follows:

\begin{enumerate}
\item For every $i = 1, . . . , r$ take a complete finite-volume non-compact
pinched negatively curved $n_i$-manifold $V_i$, where $2 \leq n_i \leq n$. 

\item Denote by $M_i$ the compact, smooth manifold with boundary obtained by removing from $V_i$ a {(nonmaximal)} horospherical open neighborhood of each cusp, thereby \emph{truncating} it. 

\item Take fiber bundles $Z_i\to M_i$, which are trivial in a neighborhood of $\partial
M_i$ and with fiber an infranilmanifold $N_i$ of dimension $n-n_i$, i.e. $N_i$ is
diffeomorphic to $\widetilde{N_i}/\Gamma_i$, where $\widetilde{N_i}$ is a simply connected
nilpotent Lie group, and $\Gamma_i$ is an extension of a lattice $L_i\subset
\widetilde{N_i}$ by a finite group. So the finite cover $N'_i$ of $N_i$ with covering
group $H_i=\Gamma_i / L_i$ is a nilmanifold $\widetilde{N_i} / L_i$. Assume that the
structure group of the bundle $Z_i\to M_i$ reduces to a subgroup of affine transformations
of  $N_i$.

\item Fix a complete pairing of diffeomorphic boundary components between
distinct $Z_i$'s, provided one exists, and glue the paired boundary components
using diffeomorphisms smoothly isotopic to affine diffeomorphisms of the boundaries, to obtain a connected manifold of dimension $n$.
\end{enumerate}

We will call the $Z_i$'s the {\em pieces} of $M$, and whenever $\dim(M_i)=n$, we say $Z_i=M_i$ is a {\em pure} piece.
\end{Definition}

\begin{Example}\label{ex::FLS}
The manifold $M$ constructed in Example \ref{ex::geoms} belongs to the family of manifolds investigated by Frigerio-Lafont-Sisto \cite{FLS15}. 
\end{Example}

\begin{Example}\label{ex::HGM}
 With $M_1$ and $M_2$ as in Example \ref{ex::geoms}, now define $Z_3$ as a fiber bundle $T^2\to Z_3\to M_1$, and $Z_4$ as a fiber bundle $S^1\to Z_4\to M_2$, both trivial in boundary collar neighborhoods. Glue $Z_3$ to $Z_4$ using an affine transformation $A:T^3\to T^3$, extended to---smaller---collar neighborhoods of their boundaries. These manifolds fall into the family described in Definition \ref{def::HGM}. Observe that they will only be part of the family described by Frigerio-Lafont-Sisto \cite{FLS15} when $A$ is the identity map, and the fibration structures are trivial products. 
\end{Example}

The JSJ decomposition in dimension three also has an analog for higher dimensional nonpositively curved manifolds. For example, consider a smooth orientable compact and connected  $n$-manifold $M$ with nonpositive curvature and convex boundary. Leeb and Scott showed that $M$ admits a geometric decomposition analogous to the topological Jaco-Shalen-Johanson \cite{JS,J} torus decomposition. Moreover, they proved that either $M$ admits a flat metric or $M$ decomposes along totally geodesic codimension one submanifolds, flat in the metric induced from $M$. The resulting pieces of this decomposition, termed the \emph{geometric characteristic splitting}, are either Seifert fibered or codimension-one atoroidal.

\begin{Definition}\label{def::Seifert-LS}\cite{LS}{ A manifold $N$ of dimension $n$ is {\bf Seifert fibered} if $N$ is a Seifert bundle over a 2--dimensional orbifold with fiber a flat $(n-2)$-manifold.

So if $N$ is Seifert fibered, then it is foliated by $(n-2)$--dimensional closed flat manifolds, each leaf $F$ has a foliated neighborhood $U$ which has a finite cover whose induced foliation is a product $F \times D^2$.

A manifold $N$ of dimension $n$ is {\bf codimension-one atoroidal} if any $\pi_1$--injective map of an $(n-1)$--torus into $M$ is homotopic into the boundary of $N$.}
\end{Definition}

With this terminology, we can prove:

\begin{Theorem}\label{thm::noEinst-hgms-NPC}

\begin{enumerate}

\item Let $X$ be a compact closed higher graph $4$--manifold without pure pieces and with a proper graph decomposition into at least two pieces. Then $X$ does not admit an Einstein metric.

\item Let $M$ be a closed orientable smooth nonpositively curved $4$--manifold. Assume that the Leeb-Scott geometric characteristic splitting of $M$ is non-trivial and consists of Seifert fibered pieces (according to Definition \ref{def::Seifert-LS}) and of pieces whose fundamental group have non-trivial center. Then $M$ does not admit an Einstein metric. 

\end{enumerate}

\end{Theorem}

Our results depend on the existence of sequences of collapsing metrics with bounded curvature, either explicitly constructed or guaranteed by polarized $\Fs$--structures and the work of Cheeger and Gromov \cite{CG}. These polarizations imply that the relevant terms of the Hitchin-Thorpe inequality vanish, thus saturating it. Therefore by the classification theorem of Hitchin \cite{hit}, we can exclude manifolds whose fundamental group has excessive growth. 

Recently, related results were presented by Di Cerbo \cite{D21} for the family of {\it extended} graph $4$--manifolds of Frigerio-Lafont-Sisto. He showed that all extended graph $4$--manifolds do not admit Einstein metrics. As we explained in Example \ref{ex::HGM} above, there exist higher graph $4$--manifolds that do not fit into this family. The results we prove apply to a larger family of graph manifolds, beyond the extended graph manifolds of Frigerio-Lafont-Sisto that were the focus of Di Cerbo's paper. Moreover, our techniques for proving these results are somewhat different. Our central point is that the same mechanism and arguments cover all of the families of manifolds included here. We use collapsing with bounded curvature, or in other cases the polarization of the $\Fs$--structures, to conclude that the Euler characteristic and signature vanish, using Hirzebruch's formula involving the first Pontryagin number. Then the computations of growth of fundamental groups are carried out in each relevant case. The existence of real hyperbolic metrics is ruled out by Lemma \ref{lem::no-hyp}. A similar strategy was already hinted at by LeBrun \cite{leb-2009}. Furthermore, as a referee pointed out, the calculation of Euler characteristics and signatures of geometric 4-manifolds by Johnson and Kotschick is relevant \cite{JoKo}. 

Here, we have systematically exploited this method, in combination with our previous results, without always relying on assumptions about complex or symplectic structures. In this way, we produce a list of manifolds that do not admit Einstein metrics that have not appeared elsewhere.  We thus address a basic issue, mentioned in Anderson's survey of Einstein $4$--manifolds \cite{Anderson}; the extent to which the Thurston view of $3$--manifolds carries over to dimension $4$, by showing that several families of collapsing manifolds do not admit Einstein structures. 

Finally, even though some of these results may be known to some experts, we hope to contribute with this note a panoramic view that surveys interesting families of manifolds that do not admit Einstein metrics using the same methods.

It is not clear presently how to deal with the combination of complex hyperbolic pieces in higher graph $4$-manifolds. This would cover all of the possible higher graph manifold pieces in dimension four. We consider this case to be outside of the scope of this note (these are non-collapsing pieces). 

In section \ref{sec::2}, we include the necessary results on collapsing and growth of fundamental groups, and show necessary lemmata. Section \ref{sec::3} contains the proofs of our main results. 

\section{Ancillaries \& Preliminaries}\label{sec::2}

\subsection{Collapsing}

A manifold $M$ is said to collapse with curvature bounded from below if it supports a sequence of metrics $g_i$ with uniformly bounded sectional curvature and volumes converging to zero as $i\to \infty$. 
An $\Fs$--structure on a closed manifold $M$ is defined by the following:
\begin{enumerate}
\item{ A finite open cover $\{ U_1, ..., U_{N} \} $;}
\item{ $\pi_{i}\co\widetilde{U_{i}}\rightarrow U_{i}$ a finite Galois covering with a group of deck transformations $\Gamma_{i}$, $1\leq i \leq N$;}

\item{ A smooth torus action with finite kernel of the $k_{i}$-dimensional torus, \\ $\phi_{i}\co T^{k_{i}}\rightarrow {\rm{Diff}}(\widetilde{U_{i}})$, $1\leq i \leq N$;}

\item{ A homomorphism $\Psi_{i}\co \Gamma_{i}\rightarrow {\rm{Aut}}(T^{k_{i}})$ such that
\[ \gamma(\phi_{i}(t)(x))=\phi_{i}(\Psi_{i}(\gamma)(t))(\gamma x) \]
for all $\gamma \in \Gamma_{i}$, $t \in T^{k_{i}}$ and $x \in \widetilde{U_{i}}$; }

\item{ For any finite sub-collection $\{ U_{i_{1}}, ..., U_{i_{l}} \} $ such that  $U_{i_{1}\ldots i_{l}}:=U_{i_{1}}\cap \ldots \cap U_{i_{l}}\neq\emptyset$ the following compatibility condition holds: let $\widetilde{U}_{i_{1}\ldots i_{l}}$ be the set of points $(x_{i_{1}}, \ldots , x_{i_{l}})\in \widetilde{U}_{i_{1}}\times \ldots \times \widetilde{U}_{i_{l}}$ such that $\pi_{i_{1}}(x_{i_{1}})=\ldots = \pi_{i_{l}}(x_{i_{l}})$. The set $\widetilde{U}_{i_{1}\ldots i_{l}}$ covers $\pi_{i_{j}}^{-1}(U_{i_{1}\ldots i_{l}}) \subset \widetilde{U}_{i_{j}}$ for all $1\leq j \leq l$. Then we require that $\phi_{i_{j}}$ leaves $\pi_{i_{j}}^{-1}(U_{i_{1}\ldots i_{l}})$ invariant and it lifts to an action on $\widetilde{U}_{i_{1}\ldots i_{l}}$ such that all lifted actions commute. }
\end{enumerate}

An $\Fs$--structure is {\em polarized} if the smooth torus action $\phi_{i}$ above is fixed point free for every $U_{i}$. 

Gromov defined the minimal volume MinVol$(M)$ to be the infimum of $\Vol(M, g)$
over all smooth metrics $g$ whose sectional curvature is absolutely bounded by 1. Cheeger and Gromov showed that the existence of a polarized $\Fs$--structure on a manifold $M$ implies the minimal volume ${\rm MinVol}(M)$ vanishes \cite{CG}. When ${\rm MinVol}(M)$=0, Chern-Weil theory implies that the Euler characteristic and the signature of $M$ vanish (cf. \cite{gromov(1982)}).

We will now recall the results used in the proofs of the main Theorems. They rely on constructions of sequences of collapsing metrics or polarized $\Fs$--structures.

\begin{Theorem} \label{thm::8-13}
\begin{enumerate}
    \item \label{thm::Fst-AGT}{\cite{S09}}
 Let $M$ be a closed orientable smooth $4$--manifold modeled on one of the Thurston geometries $ \mathbb{S}^{3}\times \mathbb{E},\,
 \mathbb{N}il^{3}\times \mathbb{E},\,
\mathbb{N}il^{4}, \, $\\ $ \mathbb{H}^{3}\times \mathbb{E},\, \Sl\times \mathbb{E},\, \mathbb{H}^{2}\times\mathbb{E}^{2},\,  \mathbb{S}ol^{4}_{1},\, \mathbb{S}ol^{4}_{0},$ or $\mathbb{S}ol^{4}_{m,n}.$ Then $M$ admits a polarized $\Fs$--structure.
    \item \label{prop::EllSfaces-TA}{\cite{PS10}}
Compact complex elliptic surfaces of zero Euler characteristic admit a polarized $\Fs$--structure.
    \item \label{thm::Fst-S1fols-Seifert4mfds}{\cite{S09}}
Let $M$ be a closed orientable smooth $4$--manifold that is either foliated by geodesic circles or admits a Seifert fibration. Then $M$ admits a polarized $\Fs$--structure.
    \item \label{thm::Fst-AGT2}{\cite{S09}}
 Let $M$ be a closed orientable smooth $4$--manifold which admits a geometric decomposition in the sense of Thurston into at least two pieces. If all of the geometric pieces of $M$ are modeled on geometries $\mathbb{H}^{3}\times \mathbb{E},\, \Sl\times \mathbb{E},$ or $\, \mathbb{H}^{2}\times\mathbb{E}$. Then $M$ admits a polarized $\Fs$--structure.
    \item \label{thm::colmet-HGMs-IMRN}{\cite{CS19}}
Let $X$ be a compact closed higher graph $4$--manifold without pure pieces. Then {\rm MinVol}$(X)=0$.
    \item \label{thm::Fst-GD}{\cite{S10}}
Let $M$ be a closed orientable smooth nonpositively curved $4$--manifold. Assume that the Leeb-Scott geometric characteristic splitting of $M$ is non-trivial and consists of Seifert fibered pieces and pieces whose fundamental group has  non-trivial center. Then $M$ admits a polarized $\Fs$--structure.
\end{enumerate}

\end{Theorem}

We next include a general result that will apply to all of the collapsing $4$--manifolds we consider in this paper. 

\begin{Lemma}\label{lem::chi-tau-zero} Let $M$ be a closed smooth oriented $4$--manifold that admits a sequence of smooth Riemannian volume collapsing metrics with bounded sectional curvatures. Then $\chi(M)=\tau(M)=0$.
\end{Lemma}

\begin{proof}

The existence of such a sequence of metrics implies MinVol$(M)=0$. In turn, by Chern-Weil theory (see \cite{gromov(1982)}), this implies the Euler characteristic and the Pontryagin numbers $p_i(M)$ of $M$ are all zero, so $\chi(M)=0$. Recall that by Hirzebruch's signature theorem $\tau(M)=\frac{1}{3}\cdot p_1(M)$ \cite{Hirz}. Therefore $\tau(M)=0$.
\end{proof}

\subsection{Simplicial Volume}\label{sec::simpvol}

Gromov defined the simplicial volume of a manifold as the infimum of $|r_i|$ with $r_i$ the coefficients of a real cycle that represents the fundamental class of $M$ \cite{gromov(1982)}. A famous result due to Thurston and Gromov states that the simplicial volume of a hyperbolic manifold is proportional to its hyperbolic volume, and in particular, it is positive. 

On the one hand, a hyperbolic structure on a manifold gives it positive simplicial volume. But, on the other hand, manifolds with an $\Fs$--structure have zero simplicial volume \cite{gromov(1982), PP}. Therefore, none of the manifolds included in Theorem \ref{thm::8-13} in the previous subsection can admit a real hyperbolic structure. We state this fact explicitly because we use it in the proofs below. 

\begin{Lemma}\label{lem::no-hyp}
Let $M$ be a manifold featured in Theorem \ref{thm::8-13}. Then $M$ does not admit a hyperbolic metric.
\end{Lemma}

\subsection{Group growth}

 A function $a: [0, \infty ) \to  {\bf R}$ is called a growth function if $a(0) \geq 1$, it is monotonically
increasing, and $a$ is submultiplicative. For $\alpha > 0 $, the growth function $(r+1)^{\alpha}$ has {\em polynomial} growth, and the growth function $e^{\alpha r}$ has {\em exponential} growth.

\subsubsection{Growth of generating sets}
Let $S$ be a finite and symmetric generating set of a group $\Gamma$. Write 
$N_{S}(m)$ for the number of elements of $\Gamma$ that can be expressed as a word of
length at most $m \in {\bf N }\setminus \{ 0 \}$ in S. Then $N_{S}$ is monotonically increasing, $N_{S}(0) = 1$ and
$N_{S}(m + n) \leq N_{S}(m)N_{S}(n)$.

 Setting $C_{a}=N_{S}(1)$ makes $N_{S}(\lfloor r \rfloor)$ a growth function, because 
 \[
 \lfloor r + s \rfloor \geq \lfloor r \rfloor + \lfloor s \rfloor + 1.
 \]
 
 \subsubsection{Growth types of groups}
 
 The growth {\em type} of a group is independent of the choice of generating set $S$. For example, if the growth of a particular generating set is exponential, then the growth of any other generating set will also be exponential. Therefore the concepts of growths of exponential or polynomial type are well defined for groups. 

\begin{Example} Consider the free group with two generators $\mathbb{F}_2$ with generating set $S=\{a,b\}$. The number of elements in $\mathbb{F}_2$ of length at most zero is one, i.e., $N_S(0)=1$. The number of elements in $\mathbb{F}_2$ of length at most one is five. Using this, we can show that the number of elements in $\mathbb{F}_2$ of length at most $n$ equals $N_S(n)=1+2(3^n+1)$. Therefore $\mathbb{F}_2$ has exponential growth. 
\end{Example}

In general, a free group with more than 2 generators has exponential growth.

Milnor first investigated the relationship between the growth of fundamental groups and curvature \cite{Milnor}. In particular, Milnor showed that negatively curved manifolds have fundamental groups of exponential growth. With respect to some of the manifolds with geometric decompositions that interest us here, Paternain and Petean \cite[Lemma 5.2]{PP} showed:

\begin{Lemma}\label{lem::expgwth-geomT}{\cite{PP}}
Let $M$ be a closed orientable smooth $4$--manifold modeled on one of the Thurston geometries $ \mathbb{H}^{3}\times \mathbb{E},\, \Sl\times \mathbb{E},\, \mathbb{H}^{2}\times\mathbb{E}^{2},\,  \mathbb{S}ol^{4}_{1},\, \mathbb{S}ol^{4}_{0},$ or $\mathbb{S}ol^{4}_{m,n}.$ Then $\pi_1(M)$ has exponential growth.
\end{Lemma}

The next result provides a specific bound for the exponential growth rate of a group. Define the exponential growth rate as \[\omega(S)=\lim_{t\rightarrow}\sqrt[k]{(N_S(k))}.\]

\begin{Theorem}\label{thm::expgrowth}{\cite{HarpeB}}
	Let $A,B$ be two finitely generated groups and let $C$ be a subgroup of both $A$ and $B$. Assume that $\left( \left[A:C\right]-1\right)\left( \left[B:C\right]-1\right)\geq 2$. Then $\omega(A\ast_C B)\geq \sqrt[4](2)$. In particular, the free product with amalgamation $A\ast_C B$ is of uniformly exponential growth. 
\end{Theorem}

The proof uses the fact that $A\ast_C B$ acts on a bipartite tree $X$.  An application of the ping-pong lemma between the vertex set made up of the disjoint union of the coset spaces $A\ast_C B/A$ and $A\ast_C B/B$ yields the desired growth. Observe that Theorem \ref{thm::expgrowth} is also valid for all HNN-extensions of a group.

Let $X$ be a compact closed higher graph manifold as in Definition \ref{def::HGM}. Suppose, to begin with, that $X$ has exactly two pieces $Z_i,Z_j$. Let $W_i\in \partial Z_i$ and $W_j\in \partial Z_j$ be a pair of diffeomorphic boundary components glued together by an affine diffeomorphism. These higher graph manifolds have a natural graph of groups associated with their fundamental groups. The groups at each of the vertices are the fundamental groups of the pieces $Z_i$ and the groups at the edges are the fundamental groups of each boundary. 

As a consequence of Eberlein's work \cite[Theorem 3.1]{Eberlein80}, each edge group $\pi_1(W_i)$ injects into the fundamental group of the piece $Z_i$. Therefore, the fundamental group of $X$ will be isomorphic to the amalgamated product $\pi_{1}(Z_i)\ast_{\pi_{1}(W_{i,j})}\pi_{1}(Z_j)$. Observe that, as we are using diffeomorphisms to glue the two pieces--up to isomorphism--this resulting amalgamated product is independent of the induced map on $\pi_{1}(W_{i,j})$.

In general, such a higher graph manifold gives rise to a graph of groups, iterating the process described for a pair of pieces over all contiguous pieces and their glued boundary components (cf. \cite{FLS15, CS19}).

\begin{Lemma}\label{lem::expgwth-HGMs}
Let $X$ be a compact closed higher graph $4$--manifold with a proper graph decomposition into at least two pieces. Then $\pi_1(X)$ has exponential growth.
\end{Lemma}

\begin{proof}
The proof is by induction on the number of pieces of $X$. Let $Z_i,Z_j$ be two pieces. Then, by Corollary 3.3 of Eberlein \cite{Eberlein80}, the edge group $\pi_1(W_i)$ is finitely generated, and by construction, the fundamental group $\pi_1(M_i)$ of the base of the piece $Z_i$ is also finitely generated. We then have that $\pi_1(W_i)$ is $\pi_1$--injective on $\pi_{1}(Z_i)$ and $\pi_{1}(Z_j)$. Therefore, the conditions of Theorem \ref{thm::expgrowth} are satisfied. From which it follows that the amalgamated product $\pi_{1}(Z_i)\ast_{\pi_{1}(W_{i,j})}\pi_{1}(Z_j)$ has exponential growth. Finally, as this amalgamated product is a subgroup of the fundamental group of $X$, we conclude that the fundamental group $\pi_1(X)$ also has exponential growth.
\end{proof}

Frigerio, Lafont, and Sisto mentioned that the fundamental group of an {\it extended} graph manifold has exponential growth. They divided the argument when the manifold has more than one piece and when it has a single piece.

We next include the results we will depend on for the other families of manifolds mentioned in the introduction. 

\begin{Lemma}\label{lem::expgwth-geoms}\cite[Lemma 21]{S10}
The fundamental group of a smooth $4$--manifold M with a proper geometric decomposition has exponential growth. 
\end{Lemma}

The following result is a consequence of the work of Avez \cite{Avez}, who showed that non-flat non-positively curved manifolds have fundamental groups of exponential growth.

\begin{Lemma}\label{lem::expgwth-NPC}
Let $M$ be a closed orientable smooth nonpositively curved $4$--manifold. Assume that the Leeb-Scott geometric characteristic splitting of $M$ is non-trivial. Then $\pi_1(M)$ has exponential growth.
\end{Lemma}

\section{Proofs}\label{sec::3}

First, we observe that Lemma \ref{lem::no-hyp} implies that all the manifolds mentioned in the results in the introduction do not admit a hyperbolic metric. 

We are now ready to present the proofs of the main results of this paper. 

\begin{proof}[Proof of Theorem \ref{thm::noEinst-Tgeom}]
These manifolds admit polarized $\Fs$--structures, by Theorem \ref{thm::8-13} (\ref{thm::Fst-AGT}). By Lemma \ref{lem::chi-tau-zero} their Euler characteristic and signatures are zero. Assume that $M$ does admit an Einstein metric. Hence the equality in (\ref{ht-int}) is attained. Therefore $M$ is known by the work of Hitchin to be finitely covered by a $4$-torus or a $K3$ surface \cite{hit}. Notice that the growth of the fundamental groups of these manifolds is at most polynomial. However, by Lemma \ref{lem::expgwth-geomT} for the geometries $ \mathbb{H}^{3}\times \mathbb{E},\, \Sl\times \mathbb{E},\, \mathbb{H}^{2}\times\mathbb{E}^{2},\,  \mathbb{S}ol^{4}_{1},\, \mathbb{S}ol^{4}_{0},$ or $\mathbb{S}ol^{4}_{m,n}$ the fundamental group grows exponentially. For the geometries $ \mathbb{S}^{3}\times \mathbb{E},\,
 \mathbb{N}il^{3}\times \mathbb{E},\,  $ or $
\mathbb{N}il^{4}, $ the uniqueness of geometric models shown by Hillman implies they can not be finitely covered by a $4$--torus \cite{Hill}. Moreover, they are not simply connected, and their fundamental groups are infinite. So they can not be finitely covered by a $K3$ surface either. \end{proof}

\begin{proof}[Proof of Corollary \ref{cor::infrasolv-noEins}] Compact $4$--dimensional infrasolvmanifolds with a fundamental group that is not virtually nilpotent are geometric manifolds, modeled on  $\mathbb{S}ol^{4}_{1},\, \mathbb{S}ol^{4}_{0},$ or $\mathbb{S}ol^{4}_{m,n}$, covered by Theorem \ref{thm::noEinst-Tgeom}.
\end{proof}

\begin{proof}[Proof of Theorem \ref{thm::noEinst-EllSfaces}]
Such elliptic surfaces admit polarized $\Fs$--structures by Theorem \ref{thm::8-13} (\ref{prop::EllSfaces-TA}). Therefore Lemma \ref{lem::chi-tau-zero} implies their Euler characteristic and signatures equal zero. 
The hypothesis that they are not finitely covered by either a $4$--torus or a Ricci-flat $K3$ surface means that Hitchin's theorem implies they can not admit an Einstein metric.
\end{proof}

\begin{proof}[Proof of Theorem \ref{thm::noEinst-S1fols-Seifert4mfds}]
These manifolds also admit polarized $\Fs$--structures, by Theorem \ref{thm::8-13} (\ref{thm::Fst-S1fols-Seifert4mfds}). But, by hypothesis, they do not fit into the classification given by Hitchin's theorem, so they can not admit an Einstein metric.
\end{proof}

\begin{proof}[Proof of Theorem \ref{thm::noEinst-geoms}]
These manifolds also admit polarized $\Fs$--structures, by Theorem \ref{thm::8-13} (\ref{thm::Fst-AGT2}). Again, Lemma \ref{lem::chi-tau-zero} implies their Euler characteristic and signatures all vanish. Their corresponding fundamental groups all have exponential growth for two reasons. The hyperbolic $\mathbb{H}^{3}, \mathbb{H}^{2}$, factors as well as the $\Sl$ factor force the pieces that have these geometries to have exponential growth. Moreover, the graph structure of $\pi_1(M)$ can also be described as a semi-direct product whose growth is exponential, as shown in Lemma \ref{lem::expgwth-geoms}. We conclude that they can not admit any Einstein metric.  
\end{proof}

\begin{proof}[Proof  of Theorem \ref{thm::noEinst-hgms-NPC}]
1) These manifolds have zero minimal volume as shown in Theorem \ref{thm::8-13} (\ref{thm::colmet-HGMs-IMRN}), and they thus collapse with bounded sectional curvature. By Lemma \ref{lem::chi-tau-zero} their Euler characteristic and signatures are zero. Their fundamental groups have exponential growth (see Lemma \ref{lem::expgwth-HGMs}). The same reasoning as in the previous Theorems above completes the proof. 

2) These manifolds admit polarized $\Fs$--structures, by Theorem \ref{thm::8-13} (\ref{thm::Fst-GD}). Therefore Lemma \ref{lem::chi-tau-zero} implies that their Euler characteristic and signatures are zero. Moreover, the growth of their fundamental groups is also exponential, as shown in Lemma \ref{lem::expgwth-NPC} above. As in the previous cases, these conditions impede the existence of Einstein metrics, as claimed. 
\end{proof}

We thank an anonymous reviewer for pointing out that a proof of Theorem \ref{thm::noEinst-Tgeom} also follows from  Berger's inequality \cite{Berger} and the computations of Euler characteristics worked out by Wall \cite{Wall}. We also thank another anonymous reviewer for leading us to the work of Avez \cite{Avez}, it simplified one of our proofs.

\providecommand{\bysame}{\leavevmode\hbox to3em{\hrulefill}\thinspace}
\bibliographystyle{amsplain}

\end{document}